\documentclass[12pt,reqno]{amsart}

\usepackage[colorlinks=false]{hyperref}
\usepackage{amssymb}
\usepackage{amsmath}
\usepackage{amsfonts}
\usepackage{float}
\usepackage{graphicx}
\usepackage{amsthm}
\usepackage{amscd}
\usepackage{calrsfs}
\usepackage{latexsym}
\usepackage{epsf}
\usepackage{oldlfont} 
\usepackage{xcolor} 

\graphicspath{{converted_graphics/}}

\begin{document}

\theoremstyle{plain}
\newtheorem{theorem}{Theorem}

\theoremstyle{definition}
\newtheorem{definition}[theorem]{Definition}
\newtheorem{example}[theorem]{Example}
\newtheorem{conjecture}[theorem]{Conjecture}

\theoremstyle{remark}
\newtheorem{remark}[theorem]{Remark}

\newcommand{\li}{\mathop{\mathrm{li}}}
\newcommand{\seqnum}[1]{\href{https://oeis.org/#1}{#1}}  

\def \bN{{\mathbb N}}
\def \BB{{\cal B}}
\def \TT{{\cal T}}
\def \II{{\cal I}}

\title[Large gaps between primes in arithmetic progressions]{Large gaps between primes in arithmetic progressions---an empirical approach}

\author{Martin Raab}
\begin{abstract} \noindent
An overview of the results of new exhaustive computations of gaps between primes in arithmetic progressions is presented. 
We also give new numerical results for exceptionally large least primes in arithmetic progressions.
\end{abstract}

\maketitle

\vspace{8mm}

{\small
	
}

\vspace{19mm}


\section{Past gaps and present goals}

What is the largest gap between two consecutive prime numbers below a given bound~$x$? This
rather simple question has proven to be quite intriguing and way more difficult to answer than one
might be inclined to think. The blurry picture of what is proven, what is heuristically expected and
what has been explicitly calculated left mathematicians unsatisfied for more than a century, and will
probably continue to do so for years to come.

Large gaps between primes have been sought and documented since at least the nineteenth century \cite{Glaisher}.
To date, the maximal gaps between prime numbers below $x$ have been found for all $x<2^{64}$ \cite{Nicely_gaps},
and even---pending computational verification---for some $x>2^{64}$ \cite{Loizides2021}.  
In light of the calculation time required to improve on the bound $x$, new results are 
few and far between, so it is natural to generalize the search to prime gaps in an arithmetic progression 
where it is easier to gather more data quickly. After all, a decent amount of empirical data
serves as a touchstone in comparing heuristical results and conjectures from the past and the future.

We will concentrate on the aspect of data computed with respect to gaps between primes in
arithmetic progressions, expanding on what has been done in this area, and put to the test some of the
predictions that have been made in the past. Specifically, we focus on the relation between the
maximal gaps below $x$ and $\log^2 x$.

\newpage

\section{Definitions}
\vspace{5mm}

{\small

\noindent
\begin{tabular}{ll}
\vspace{2mm}
$p_n$ & the $n$-th prime; ${p_n}\in \{2,3,5,7,11,13,17\ldots\}$ \\
\vspace{2mm}
$q$ & integer; the common difference in the arithmetic progression $r+q\cdot j$, \,$j=0,1,2,\ldots$ \\
\vspace{2mm}
$r$ & integer coprime to $q$, \,$1\leq r<q$ \\
$g$ & the gap between two primes $p^-$ and $p^+$ in the arithmetic progression $r+q\cdot j$; \\
\vspace{2mm}
& $g=p^+-p^-=q\cdot k$ \\
$p^-$ & the prime at the start of a gap of length $g$ in the arithmetic progression $r+q\cdot j$, \\
\vspace{2mm}
& also referred to as the lower bounding prime of a gap \\
$p^+$ & the prime at the end of a gap of length $g$ in the arithmetic progression $r+q\cdot j$, \\
\vspace{2mm}
& also referred to as the upper bounding prime of a gap \\
\vspace{2mm}
$p_{\min}(q,k)$ & least prime such that there is a gap of length $q\cdot k$ between $p^-$ and $p^+$ \\
\vspace{2mm}
$g_{\max}(x)$ & the maximal gap between consecutive primes $\leq x$ \\
\vspace{2mm}
$g_{\max}(x,q)$ & the maximal gap between primes $\leq x$ in an arithmetic progression \\
$m$ & the ``merit'' of a gap $g$; expected average number of primes between $p^-$ and $p^+$, \\
\vspace{2mm}
& further explained in sec.\,\ref{secPM} \\
\vspace{2mm}
CSG ratio & the Cram\'er--Shanks--Granville ratio which relates $g_{\max}(x)$ to $\log^2 x$ \\
\vspace{2mm}
$\rho(p,q)$ & a redefinition of the CSG ratio, further explained in sec.\,\ref{secPM}; $\rho(p,q)=m^2\cdot\varphi(q)/g$ \\
\vspace{2mm}
$\rho_{\max}(x,q)$ & the maximum value of $\rho(p,q)$ for all $p\leq x$ \\
\vspace{2mm}
$G(x)$ & the smooth part of Riemann's prime counting function, further explained in sec.\,\ref{secPM} \\
\vspace{2mm}
$\pi(x)$ & the prime counting function; the number of primes $p_n\leq x$ \\
\vspace{2mm}
$\varphi(q)$ & Euler's totient function \\
\vspace{2mm}
$\log x$ & the natural logarithm of $x$ \\
\vspace{2mm}
$\log^2 x$ & the same as $(\log x)^2$ \\
\vspace{2mm}
$\lfloor x\rfloor$ & the floor function; the largest integer $\leq x$ \\
\vspace{2mm}
$|x|$ & the absolute value of $x$; $|x|=-x$ for $x<0$ \\
\end{tabular}
}

\newpage

\section{Looking for the perfect measure} \label{secPM}

A key indicator here will be what is most commonly known as the Cram\'er ratio or even 
the Cram\'er--Shanks--Granville ratio. Cram\'er \cite{Cramer} conjectured, on probabilistic grounds, that
$$
\limsup\limits_{n\to\infty} \frac{p_n-p_{n-1}}{\log^2 p_n} = 1,
$$
where $p_n$ is the $n$th prime number. This is achieved by working with a random model that employs a sequence
of independent random variables $\xi_n$ where, for $n\geq3$,
$$
{\mathbb P}(\xi_n=1)=\frac{1}{\log n} \quad \mbox{ and } \quad {\mathbb P}(\xi_n=0)=1-\frac{1}{\log n}.
$$

Shanks \cite{Shanks1964} reformulated Cram\'er's conjecture as $g_{\max}(x) \sim \log^2 x$, 
with $g_{\max}(x)$ being the maximal gap between consecutive primes below $x$.

Granville \cite{Granville}, in contrast, arrived at a 
somewhat different result, using arguments of divisibility by small primes, that there should be
infinitely many instances where $g_{\max}(x) \gtrsim 2e^{-\gamma} \log^2x \approx 1.1229 \log^2x$.
The argument that the random variables $\xi_n$ in Cram\'er's model are not strictly independent
plays an essential role in what is nowadays the most widely accepted theory regarding gaps between primes.
For more details, see also Pintz \cite{Pintz2007}.

Still, the actual data collected thus far seems to be more in favor of Cram\'er and Shanks.
\vspace{8mm}

Since the same sieving rules as used in the sieve of Eratosthenes apply to both the primes in
general and primes in an arithmetic progression (AP),\footnote{
  \,Although this might not entirely be taken for granted: if the common difference in the AP is a large
  primorial, or in general has many small factors, it might have a measurable impact on the outcome. Further
  study is needed. Common differences equal to powers of 2, however, should be unbiased in this regard as they
  don't---on average---affect any residue classes modulo any odd prime.
} 
it is reasonable to suggest that, by and large,
results for primes in AP can be compared to statistical models like the ones that Granville \cite{Granville} or
Banks, Ford and Tao \cite{Banks-Ford-Tao} use for the study of prime gaps.

Throughout the paper, the ratio as described above will be measured by a function, denoted by
$\rho(p,q)$ or simply $\rho$ when $(p,q)$ is not specified, asymptotic to the above mentioned ratio. Particularly,
we look at primes in AP with common difference $q$, such that $p = r + jq$, 
where $q$ and $r$ are coprime integers, $0 < r < q$, with $q$ even, and 
$j \in {\mathbb N}_0$. A gap between two consecutive primes $p^{-}$ and $p^{+}$ 
in an AP is conventionally associated with a measure of the form\footnote{
\,cf.\ inequality (34) in \cite{Kourbatov_Wolf_2019}
}
$$
\rho^* = \frac{p^{+}-p^{-}}{\varphi(q)\log^2 p^{+}},
$$
where $\varphi(q)$ is Euler's totient function, i.e.\ the number of positive integers not exceeding $q$ 
and coprime to $q$. When $q \leq 2$, there is a different conventional measure which uses the squared logarithm of 
the lower bounding prime of the gap,
$$
\rho' = \frac{p^{+}-p^{-}}{\log^2 p^{-}}.
$$

When applied to gaps between primes in AP though, the corresponding measure, with $\varphi(q)$ in the
denominator, can produce arbitrarily large values for very small $p^{-}$ (have a try, just for fun, at $p^{-} = 3$
and $q = 6336488$). Yet, using the upper bounding prime, $\rho^*(p, q)$ is similarly 
inconclusive---although more appropriate---for very small $p^{+}$, as the respective gap is linked to
the asymptotic density of primes only at the point where the upper bounding prime $p^+$ is located.
Moreover, any such conventional measures transfer the ``flaw'' of the logarithmic integral $\li(x)$
overestimating $\pi(x)$ for most $x$ (see \cite{Johnston} for a very recent paper on this issue).
To adequately process the results of the exhaustive computation, we demand a measure that captures
the purpose of the conventional measure, the Cram\'er--Shanks--Granville (CSG) ratio $g/\log^2 p$,
while taking into account the exact average local distribution of primes in the ranges of concern. \\

Enter Gram's variant of Riemann's prime counting formula \cite{Gram}, that is, omitting the non-trivial zeroes
of the Riemann zeta function. Here we have a non-decreasing,
smooth function for $x > 1$ which approximates the prime counting function $\pi(x)$ better than any other
historically established approximation (by e.g. Legendre or Gauss) for a vast majority but especially for
small values of $x$:

\begin{equation} \label{eq1}
G(x) = 1 + \sum_{n=1}^\infty \frac{\log^n x}{n\cdot n! \cdot \zeta(n+1)}
       + \frac{1}{\pi}\arctan \frac{\pi}{\log x} - \frac{1}{\log x}.
\end{equation}\\

\vspace{-30mm}
\begin{figure}[H]
  \centering
  \includegraphics[bb=-100 0 720 468,width=7in,height=4.5in,keepaspectratio]{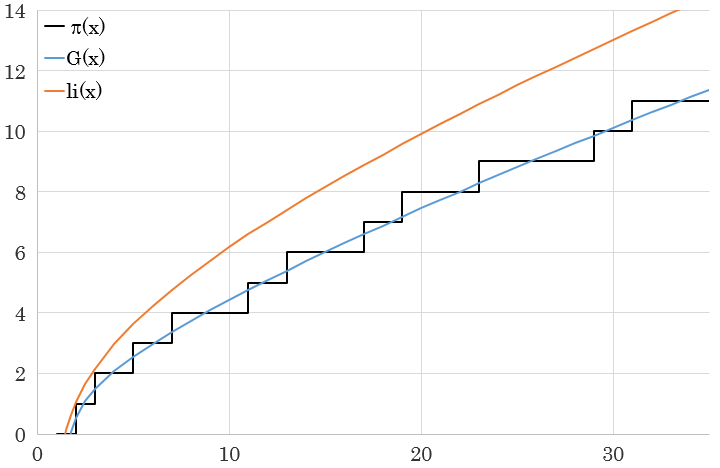}
  \vspace{3mm}
  \caption{Comparison: $\pi(x)$ vs. $G(x)$. As we go to larger $x$ values, $x=97$ is the first prime where
  $|G(x-\varepsilon)-\pi(x-\varepsilon)|>1$ for $\varepsilon<0.13$. The logarithmic integral $\li(x)$, implicitly
  used for the conventional measures, clearly is off on a small scale.
}
\label{fig1}
\end{figure}

The derivative of this function, $G'(x)$, may serve---in the spirit of Cram\'er---as the ``probability'' for
a randomly chosen integer $x$ to be prime. (Of course, $x$ is either prime or not, but practically speaking,
we assume the probability prior to choosing the exact value of $x$.) \\

The {\em merit} $m$ of a gap indicates how many primes are expected on average for a particular gap:
\begin{equation} \label{eq2}
m = \frac{q}{\varphi(q)} \sum_{j=1}^k G'(p+jq),
\end{equation}\\
where $p=p^-$ is the prime preceding the gap, and
$p+kq=p+g=p^+$ is the prime following the gap. \\

The formula for the CSG ratio can then be redefined (eliminating explicit dependence on~$p$)
as follows:
\begin{equation} \label{eq3}
\rho = m^2 \,\frac{\varphi(q)}{g},
\end{equation}\\
which is asymptotically equivalent to the conventional measures for $p\to\infty$.
Furthermore noting that, for odd $q$,
only one respective odd value of $k$ presents a gap for $p^-= 2$, this measure is only coherent if we solely
use even values of $q$; odd values of $q$ may be used instead whenever $q$ is congruent to 2 mod 4, but
then we would ``cheat'' on $\rho$ being a bit larger than otherwise. Only now it is possible to start analyzing
the data of first occurrences or maximal gaps for any $q$ we choose without it being distorted for small
values of $p$ or large values of $q$, a problem we would face with any of the commonly known
conventional formulas for the CSG ratio.

\newpage

\section{A plethora of data}

For many pairs $(q,k)$, with $k\in{\mathbb N}$ and $q$ even, we have computed the least prime 
$p_{\min}=p_{\min}(q,k)$ such that there is a gap of length $k \cdot q$ between $p_{\min}$ 
and the next prime in the AP with difference $q$. Table \ref{table1} gives some results of this computation 
for small values $k$ and $q$.

\begin{table}[H]
\caption{Least primes $p_{\min}(q,k)$ such that there is a gap of length $k \cdot q$ 
 before the next prime in the AP with difference $q$.
 For any value of $q$, it is difficult to extend this table horizontally (say, to $k > 500$); 
 however, it's very easy to compute more rows for small $k$.
\label{table1}
}
\centering
\small {
\begin{tabular}{lccccccccccc}
\hline $\hspace{1cm}{\vphantom{\displaystyle\int}}$
     & $k=1$  & $k=2$ & $k=3$ & $k=4$ & $k=5$ & $k=6$ & $k=7$ & $k=8$ & $k=9$ & $k=10$ & $k=100$  \\
[0.5ex]\hline
$q=2$  \vphantom{\fbox{$1^1$}} 
            & 3 & 7 & 23 &  89 & 139 & 199 & 113 & 1831 & 523 & 887 & 378043979 \\
$q=4$   & 3 & 5 & 17 & 73 & 83 & 113 & 691 & 197 & 383 & 1321 & 107345389 \\
$q=6$   & 5 & 19& 43 &  283 & 197 & 521 & 1109 & 2389 & 1327 & 4363 & 680676109 \\
$q=8$   & 3 & 7 & 17 &  41 & 61 & 311 & 137 & 451 & 647 & 1913 & 63977327 \\
$q=10$ & 3 & 11 & 29 & 313 & 113 & 397 & 331 & 269 & 997 & 1129 & 881451157 \\
$q=12$ & 5 & 13 & 53 & 109 & 379 & 1109 & 457 & 2111 & 2711 & 1667 & 1906215407 \\
$q=14$ & 3 & 13 & 11 & 101 & 127 & 233 & 761 & 661 & 1091 & 1619 & 467186417 \\
$q=16$ & 3 & 5 & 19 & 127 & 17 & 137 & 139 & 449 & 617 & 1063 & 185110507 \\
$q=18$ & 5 & 7 & 59 & 107 & 269 & 631 & 727 & 677 & 1709 & 5167 & 2743266193 \\
$q=20$ & 3 & 7 & 29 & 101 & 283 & 239 & 569 & 433 & 1823 & 3257 & 618986273 \\
\hline
\end{tabular}
}
\end{table}

In table \ref{table1}, all residue classes coprime to $q$ are present. One could go further and analyze each
residue class separately. It should be noted that since the residue classes that are squares mod $q$ have
a bias to contain, under certain circumstances, less primes than those that are non-squares mod $q$ (for
more details see \cite{Granville-Martin}), it might make sense to check whether this effect is measurable in the table of
gaps. If so, the residue classes that comprise first occurrence gaps should be more numerous if said
residue classes are squares mod $q$.

For $k \leq 100$, the number of primes $p_{\min}(q,k)$ that are squares modulo $q$ are as shown 
below in table \ref{table2}:

\begin{table}[H]
\caption{The number of primes $p_{\min}(q,k)$ that are squares mod $q$, for $k \leq 100$.
\label{table2}
}
\centering
\small {
\begin{tabular}{cccccccccc}
\hline $p_{\min}\equiv$ square mod $q$? ${\vphantom{\displaystyle\int}}$
      & $q=4$ & $q=6$ & $q=8$ & $q=10$ & $q=12$ & $q=14$ & $q=16$ & $q=18$ & $q=20$  \\
\hline
\ yes  
${\vphantom{1^{1^1}}}$
      & 53 & 54 & 34 & 45 & 28 & 49 & 28 & 41 & 22 \\
no   & 47 & 46 & 66 & 55 & 72 & 51 & 72 & 59 & 78 \\
\hline
\end{tabular}
}
\end{table}

For $q = 8$, 12, 16, and 20, only 25\% of residue classes coprime to $q$ are squares, so instead of 50:50,
we might expect a 25:75 ratio here if evenly distributed. 
The impact of squares vs. non-squares mod $q$
might be subject to further study, but can be considered ancillary for the motif of this paper.

\begin{table}[H]
\caption{The maximum values $\rho_{\max}(x,q)$ up to a given point $x$ for each $q$.
The last lines show averages and weighted averages by a factor of $\varphi(q)/q$, for $q \leq 5000$.
\label{table3}
}
\centering
\small {
\begin{tabular}{lccccccccc}
\hline ${\vphantom{\displaystyle\int}} \rho_{\max}(x,q) $  & 
$x=2^{21}$& $x=2^{24}$& $x=2^{27}$& $x=2^{30}$& $x=2^{33}$& $x=2^{36}$& $x=2^{39}$& $x=2^{42}$& $x=2^{45}$  \\
\hline
$q=2$ ${\vphantom{1^{1^1}}}$
            & 0.7020 & 0.7020 & 0.7393 & 0.7393 & 0.7393 & 0.7953 & 0.7953 & 0.7975 & 0.8178 \\
$q=4$   & 0.6792 & 0.6792 & 0.7511 & 0.7511 & 0.7678 & 0.7678 & 0.8057 & 0.8156 & 0.8157 \\
$q=6$   & 0.7360 & 0.7360 & 0.7540 & 0.7540 & 0.7540 & 0.7540 & 0.7540 & 0.8043 & 0.8286 \\
$q=8$   & 0.6397 & 0.6398 & 0.7248 & 0.7248 & 0.7437 & 0.7532 & 0.7876 & 0.7876 & 0.7876 \\
$q=10$ & 0.7132 & 0.7132 & 0.7132 & 0.7677 & 0.8043 & 0.8043 & 0.8930 & 0.8930 & 0.9050 \\
$q=12$ & 0.6841 & 0.6957 & 0.7048 & 0.7371 & 0.7586 & 0.8053 & 0.8455 & 0.8455 & 0.8455 \\
$q=14$ & 0.6623 & 0.6623 & 0.6671 & 0.7694 & 0.8096 & 0.8346 & 0.8346 & 0.8346 & 0.8346 \\
$q=16$ & 0.7177 & 0.7313 & 0.8110 & 0.8110 & 0.9046 & 0.9046 & 0.9046 & 0.9046 & 0.9046 \\
$q=18$ & 0.7870 & 0.7870 & 0.7870 & 0.7870 & 0.7870 & 0.7870 & 0.7870 & 0.7870 & 0.7879 \\
$q=20$ & 0.7222 & 0.7941 & 0.7941 & 0.7941 & 0.7941 & 0.7941 & 0.7941 & 0.7941 & 0.8238 \\
\hline
\multicolumn{2}{c}{Average, $q \leq 5000 {\vphantom{1^{1^1}}}$:} 
    & 0.7310 & 0.7573 & 0.7795 & 0.7966 & 0.8099 & 0.8209 & 0.8315 & 0.8406 \\
\multicolumn{2}{c}{Weighted average:}     
    & 0.7316 & 0.7580 & 0.7799 & 0.7971 & 0.8105 & 0.8215 & 0.8321 & 0.8412 \\
\hline
\end{tabular}
}
\end{table}

In table \ref{table3} we track the maximum values $\rho_{\max}(x,q)$  of the redefined CSG ratio $\rho$  up to a given point
$x$ for each $q$.\footnote{\,By the time $\rho_{\max}(x,q)>0.8$, most of these values agree to at least three decimal
places with the conventional measure $\rho^*$ mentioned in sec.\,\ref{secPM}.}
It also makes sense to consider the weighted averages of $\rho_{\max}(x,q)$  by a factor of $\varphi(q)/q$,
since $\rho_{\max}(x,q)$ tends to be smaller on average whenever $\varphi(q)/q$ is small.

How much smaller? For a more nuanced picture, we split the data into three groups for which we obtain the respective average
values of $\rho_{\max}(x,q)$:
$$
a)\,1/2.1 < \varphi(q)/q \leq 1/2; \quad \quad b)\,1/3 \leq \varphi(q)/q < 1/2.1; \quad \quad c)\,\varphi(q)/q < 1/3.
$$
These groups are of roughly comparable size: in the long run, 33.1\% of even $q$ fall into group a), 32.8\% into group b),
and 34.1\% into group c). We're only looking at $q \leq 5000$ though, where group b) is still more dominant---in terms of
sample size---in table \ref{table4} with 884 of 2500 values (35.4\%), compared to 802/2500 for group a) and 814/2500 for group c).

\begin{table}[H]
\caption{The average maximum values $\rho_{\max}(x,q)$ up to a given point $x$ for each group of values of $\varphi(q)/q$,
$q \leq 5000$, $q$ even.
\label{table4}
}
\centering
\small {
\begin{tabular}{lccccccccc}
\hline ${\vphantom{\displaystyle\int}} \varphi(q)/q $ &
$x=2^{21}$& $x=2^{24}$& $x=2^{27}$& $x=2^{30}$& $x=2^{33}$& $x=2^{36}$& $x=2^{39}$& $x=2^{42}$& $x=2^{45}$ \\
\hline
$> 1/2.1$ ${\vphantom{1^{1^1}}}$
          & 0.7023 & 0.7361 & 0.7628 & 0.7832 & 0.7993 & 0.8139 & 0.8249 & 0.8359 & 0.8449 \\
\hline
$\begin{matrix} \geq 1/3, {\vphantom{1^{1^1}}}\\< 1/2.1 \end{matrix}$
          & 0.6932 & 0.7293 & 0.7546 & 0.7769 & 0.7954 & 0.8087 & 0.8195 & 0.8303 & 0.8393 \\
\hline
$< 1/3$ ${\vphantom{1^{1^1}}}$
          & 0.6894 & 0.7277 & 0.7550 & 0.7786 & 0.7952 & 0.8074 & 0.8184 & 0.8283 & 0.8379 \\
\hline
\end{tabular}
}
\end{table}

\newpage

The weighted average ($Ø_w$) values of $\rho_{\max}(x,q)$ over $q \leq 5000$ closely satisfy the approximation
$0.963-3.8/\log x$ in the surveyed range of $x$. But this is, mildly put, a risky
assumption, as there isn't yet any heuristical indication to back up a claim like the existence
of only finitely many gaps larger than $0.963\,\varphi(q)\,\log^2 x$ (or with a factor close to that
number).\footnote{\,However, a careful inspection of the number of gaps larger than $m\log p$ for various samples
in the range $e^{24} \leq p \leq e^{600}$ indicates a density of said gaps of $e^{-m+O(m)/\log p}$---
a crucial aberration to the regular Poisson distribution with respect to $m$ according to
Cram\'er's model, hinting that even Cram\'er may have overestimated the sizes of the largest
gaps below $x$, albeit not implicitly by a constant factor.} A possible approximation based on the
Cram\'er model would be, $Ø_w \, \rho_{\max}(x,q) \approx 1-(2.5-8\varepsilon)/(\log x)^{0.8-\varepsilon}$
for some $\varepsilon < 0.02$. It is more difficult to find a similar approximate formula with leading term 1.1229.
Eventually there might be no trivial asymptotic relation reconciling both heuristical and empirical results for $\rho_{\max}(x,q)$.

\begin{figure}[H]
  \centering
  \includegraphics[bb=0 0 800 500,width=7in,height=5in,keepaspectratio]{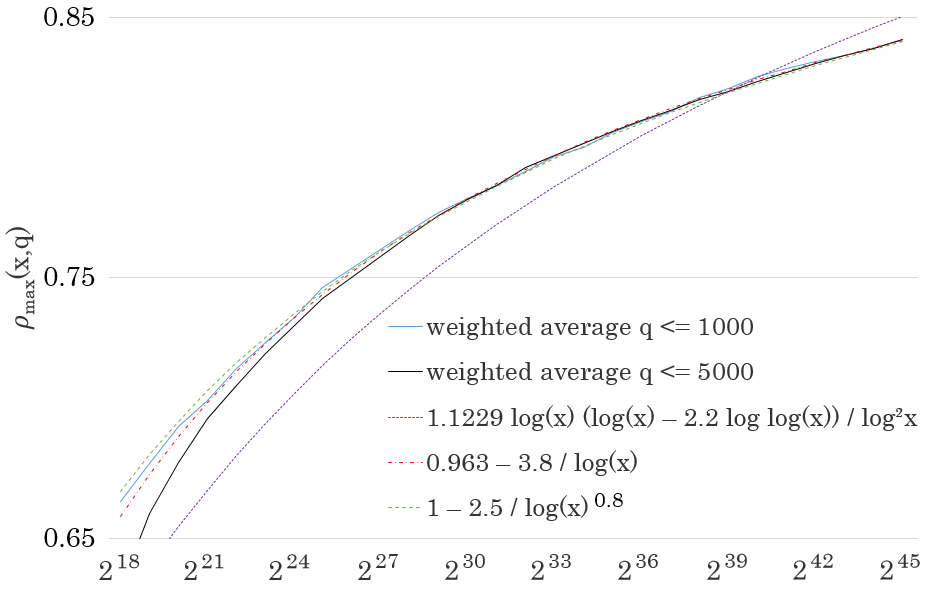}
  \caption{A formula predicting the growth of the average values of $\rho_{\max}(x,q)$ ostensibly requires
three different constants to remain in accordance with the heuristic foundations.
}
\label{fig2}
\end{figure}

Formulas of the type $a \log(x) (\log(x) - b \log \log (x))$---as suggested e.g.\ by Cadwell \cite{Cadwell} and reiterated
by Granville and Lumley \cite{Granville-Lumley}---barely hold up to the challenge of presenting a good fit for
the data thus far; see figure \ref{fig2}.

\newpage

\section{Gaps and gems}\label{secGapsGems}

\begin{figure}[H]
  \centering
  \includegraphics[bb=0 0 960 500,width=8.5in,height=4.2in,keepaspectratio]{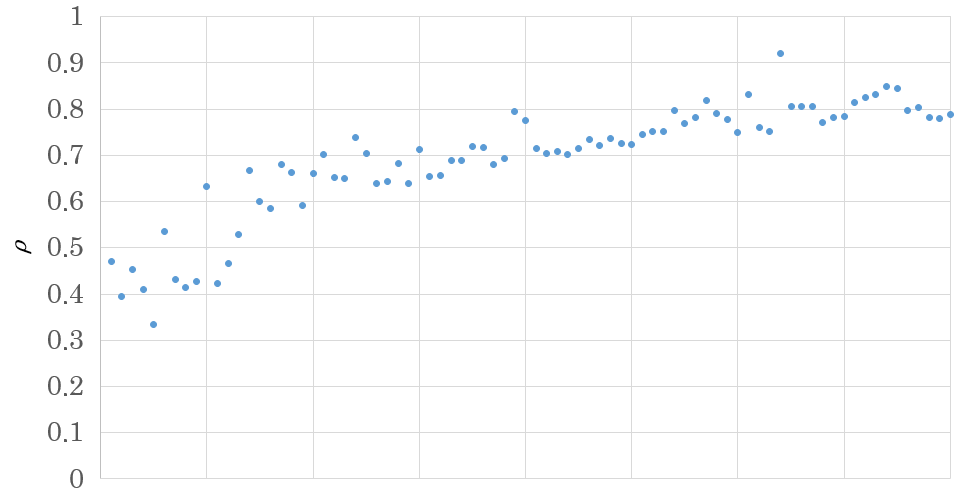}
  \caption{Graph of $\rho$ for the first 80 known maximal gaps between primes.
}
\label{fig3}
\end{figure}

Although it is conjectured that there exist infinitely many gaps between primes that are larger than
$\log^2 p_n$, computations so far leave little hope that such an exceptionally large gap will ever be found.
Nyman \cite{NymanNicely} came closest to this goal with a gap of 1132 following the prime 1693182318746371 and
$\rho \approx 0.9206$. This gap in itself was an unusually large one---seen in figure \ref{fig3} as the only point above 0.9.

Second to Nyman's gap is a gap with $\rho \approx 0.8483$ found by Oliveira e Silva {\it et al.} \cite{OliveiraeSilva2014}. 
To see how far away this is from a gap that would be on par with $\log^2(p^+)$, e Silva's gap should have had 
$p^{+}-p^{-} \ge 1700$ instead of the actual 1442. It fell short by 258, an interval in which more than six primes
of comparable size occur on average.

\vspace{-8mm}

\begin{figure}[H]
  \centering
  \includegraphics[bb=0 0 968 489,width=8.5in,height=4.2in,keepaspectratio]{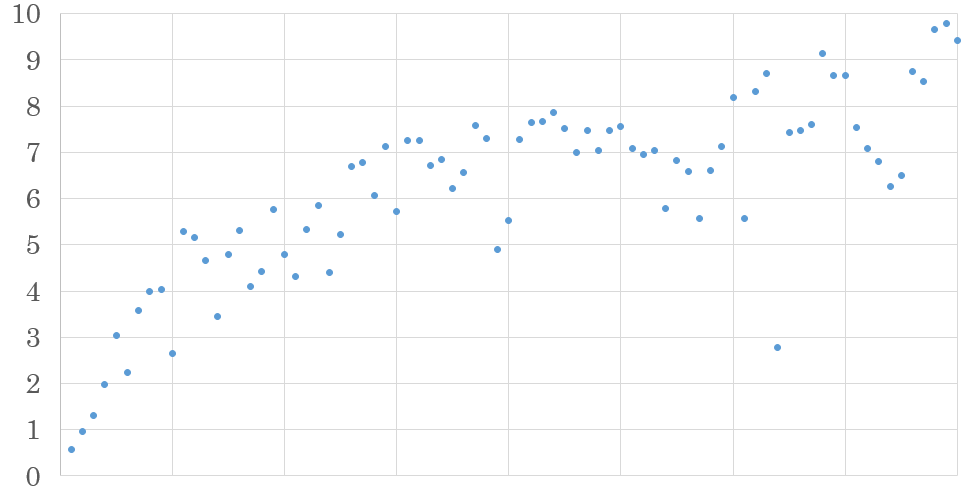}
  \caption{It would have taken an extra interval where this many more primes occur, on average, 
    to reach the magic $\log^2p$. Plotted for the first 80 known maximal gaps between primes.
}
\label{fig4}
\end{figure}

Let's have a closer look at figures \ref{fig3} and \ref{fig4}. Since new maximal gaps are already harder
to find than a needle in a haystack, it's very likely we may never see a gap between consecutive primes
where $\rho > 1$ (whether or not it exists). One could also take a guess at the growth rate in
figure \ref{fig4}---the values on the $y$-axis---and assume that maybe the sum of $1/e^y$
converges, thus likewise implying that at most finitely many gaps larger than $\log^2p$ exist.

On the other hand, the refinements of Granville's model only come into play for much larger $x$. We
may compare this to the phenomenon \cite{Engelsma} that there exist 447 numbers in an~interval
of 3159 integers coprime to arbitrarily large primorials $q\#$, hence it is possible---yet not proven---that
there are infinitely many $x$ such that $\pi(x+3159) = \pi(x)+447$ while $\pi(3159) = 446$, constituting
clusters of primes more densely packed than the primes at the start of the number line. By the 
Hardy--Littlewood conjecture \cite{Hardy-Littlewood}, the first such $x$ should be found in the vicinity 
of $10^{1198}$, which is way beyond our current computing capabilities and will likely be so for all time. \\

If we turn our attention to primes in AP, however, instances where $\rho > 1$ can indeed be found. The
first such example was discovered by Kourbatov \cite{Kourbatov_conj77} and by 2019, Kourbatov and Wolf  
\cite{Kourbatov_Wolf_2019} have listed 35 of these gems. Upon this, we continued the search for these 
exceptionally large gaps and found several hundred more. Still, they are scarce enough to keep track of 
each and every one of them, as their distribution raises more questions on the way up the search bounds.

\begin{figure}[H]
  \centering
  \includegraphics[bb=0 0 958 612,width=8.5in,height=6in,keepaspectratio]{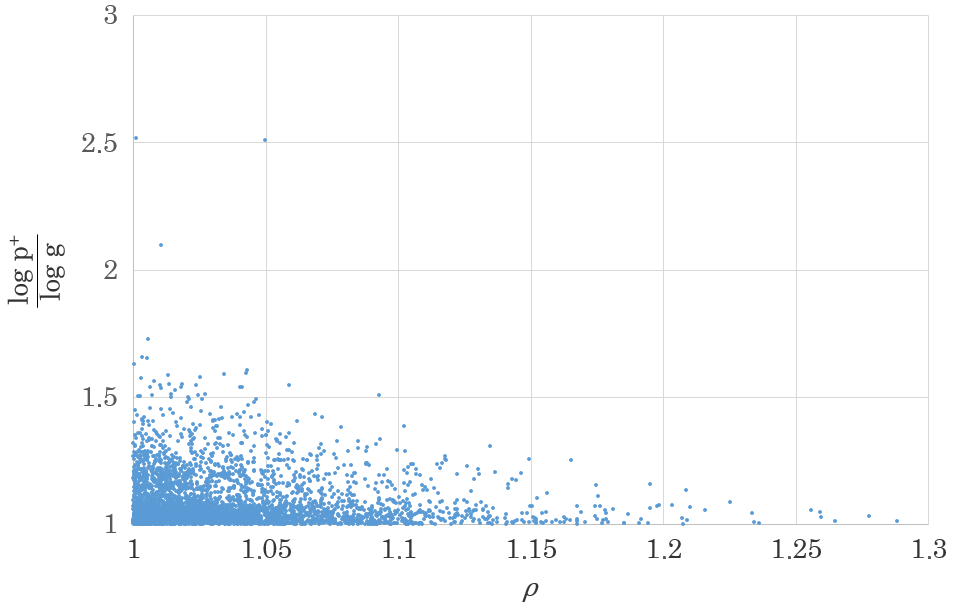}
  \caption{The ratio $\log p^+/\log g$ vs. $\rho$ for all known exceptional gaps.
}
\label{fig5}
\end{figure}

\begin{figure}[H]
  \centering
  \includegraphics[bb=0 0 958 612,width=8.5in,height=6in,keepaspectratio]{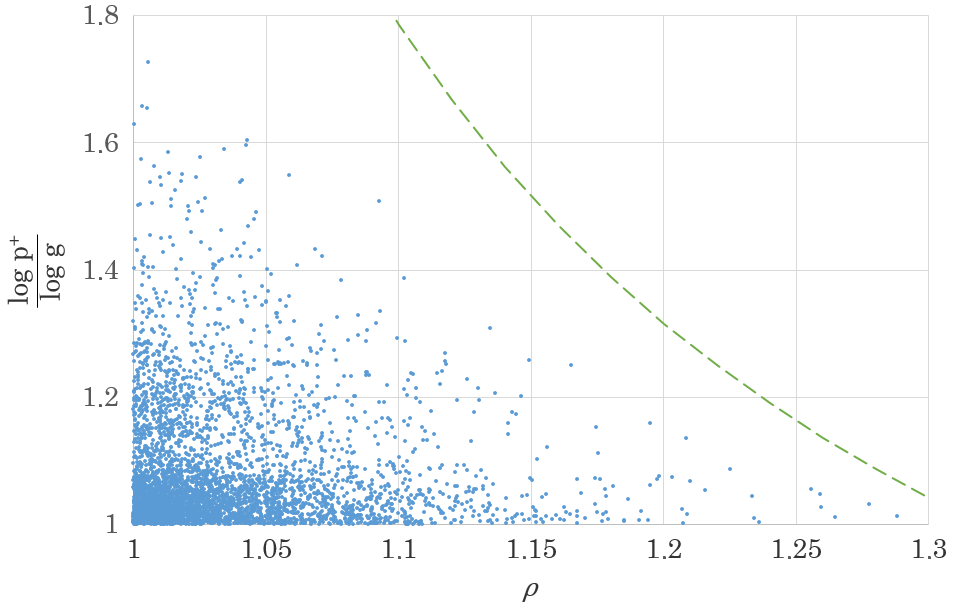}
  \caption{The correlation becomes more apparent when zooming in and concentrating on $\log p^+/\log g < 1.8$. 
The data is the same as in fig.\,5, only without the extreme examples at $q$ = \{152, 762, and 3796\},
which are also mentioned below. The dashed green curve is a speculative bound; it would certainly be
interesting to know whether arbitrarily many points come closer to or even exceed this bound or
whether the picture remains in about the current state and further exceptional gaps tend to appear
ever closer to $(1,1)$ in the bottom left corner.
}
\label{fig6}
\end{figure}

\newpage

One of the more compelling questions concerns the relation between $p$ and $q$ with respect to the
magnitude of $\rho$. For small $p^{-}$, prime factors that contribute to the gap (as they are divisors of the
composite numbers $p^{-}+jq$ for $0 < j < k$) may only just appear as the numbers within the gap get larger.
In other words, the conditions for a number being prime at the start of a gap can be different than at
the end of the gap, in the precise sense that $\lfloor\sqrt{p^-} \rfloor < \lfloor\sqrt{p^+} \rfloor$.
The change of these conditions within the gap, which are more prominent when $\log p^+$ is closer to $\log g$,
is the main reason why it is possible to find gaps with  $\rho > 1$ in the first place (anyhow considering what
is computationally feasible today).
For exceptionally large gaps with $\rho > 1$ found thus far, $\log p^+/\log g$ shows a distribution inversely
proportional to $\rho(p,q)$, as to be seen in figures \ref{fig5} and \ref{fig6} as well as in table \ref{table5}. \\

For $q$ fixed, the larger $p$ gets the more difficult it is to find gaps with $\rho$ exceeding a certain fixed
bound. Once an exceptional gap is found for a certain $q$, it is rare to find another exceptional gap for
the same $q$, only three cases are known so far:
$$
\ \ q = 28388, \ \ \quad p^- = 366870073 \ \mbox{ and } \ p^- = 5088100651;
$$
$$
\ \ q = 389104, \ \ \ \,  p^- = 461954737 \ \mbox{ and } \ p^- = 2176128499;
$$
$$
\ \ q = 1238684, \ \  p^- = 87686639 \ \mbox{ and } \ p^- = 19190651717.
$$

All these gaps have $\rho < 1.05$. \\

Even harder to find are exceptional gaps for which $p$ is larger than the square of the gap--- or, closer to the
issue of potential prime factors,  $\lfloor\sqrt{p^-} \rfloor = \lfloor\sqrt{p^+} \rfloor$ --- where $\log p^+/\log g>2$.
Only two such instances have been found so far:\footnote{\,Note also that $762 \equiv$ 2 (mod 4):
if we set $q = 381$ here, both of these gaps are also larger than $q^2$.}
$$
q = 152, \quad k = 357, \quad \ \ \ p^- = 825353008489;
$$
$$
q = 762, \quad k = 340, \quad p^- = 38943534114929.
$$

These are depicted by the data points just above 2.5 on the left in figure \ref{fig5}. \\

Another example that barely misses the condition $\lfloor\sqrt{p^-} \rfloor = \lfloor\sqrt{p^+} \rfloor$, 
yet $p$ is larger than the square of the gap, is
$$
q = 3796, \quad k = 411, \quad p^- = 9585778010467.
$$ \\

Out of all 4070 exceptional gaps found to date, of which 2257 satisfy $g > \varphi(q) \log^2 p^+$ by the
conventional criterion, 260 have $\rho > 1.1$, and a mere 17 have $\rho > 1.2$.
Statistically, there is a rough relation such that the number of gaps $g_\delta$ with $\rho > 1+\delta$ is
approximately\footnote{\,And here we likewise have a correlation between the constant in the exponent
and the ratio $\log p^+/\log g$: if said ratio is smaller than 1.1, the above relation is closer to $g_0\,e^{-27 \delta}$,
whereas for $\log p^+/\log g > 1.1$, we have $g_\delta\,\approx\,g_0\,e^{-31 \delta}$.}
$g_0\,e^{-28 \delta}$, which in practice holds well for $\delta < 0.17$, but no more so for larger $\delta$.
And for those gaps with a large $\rho$, the ratio $\log p^+ / \log g$ is more likely to be very close to 1; see table \ref{table5}.

\newpage

\begin{table}[H]
\caption{Largest value of $ \log p^+ / \log g$ for gaps with at least the given magnitude of $\rho$.
The data suggests that the largest $\rho$ can be found when $p$ is not too large compared to the gap size
$g$. These results are partially influenced by higher search bounds of $p$ for smaller $q$, though the scope
of this effect is somewhat unclear.
\label{table5}
}
\centering
\begin{tabular}{cccc}
\hline $\ \rho > {\vphantom{\displaystyle\int}}$ & $\max \frac{\log p^+}{\log g}$ & for $p^-$ = & and $q$ = \\
\hline
${\vphantom{1^{1^1}}}$
1.28 \ & 1.0123 & 209348411 & 2830474 \\
1.27 & 1.0318 & 938688203 & 4200826 \\
1.25 & 1.0554 & 1961096147 & 3128278 \\
1.20 & 1.1354 & 5120249753 & 1251242 \\
1.15 & 1.2499 & 25319877559 & 875600 \\
1.10 & 1.3871 & 48179541911 & 245466 \\
1.05 & 1.5483 & 74651827093 & 33554 \\
1.04 & 2.5105 & 38943534114929 & 762 \\
1.00 & 2.5170 & 825353008489 & 152 \\
\hline
\end{tabular}
\end{table}

The first line in table \ref{table5} has in fact $\rho = 1.28829...$. Can $\rho$ get any larger than that? We don't know.
It might be that this is, by the measure in question, the largest gap there is for any choice of $q$. It is conceivable that
$\rho(p,q)$ might exceed 1.3 for some set of pairs $(p,q)$, yet all data collected so far indicates that it won't get much
larger. A theoretical approach that might shed some light on this problem would be greatly welcomed.
\vspace{5mm}

The search depth $p$ for various levels of $q$, as of March 2023, are shown in table \ref{table6}.

\begin{table}[H]
\caption{Search bounds. Only even $q$ are examined.
\label{table6}
}
\centering
\begin{tabular}{cc}
\hline $\ q \leq {\vphantom{\displaystyle\int}}$ & $p \leq $ \\
\hline
${\vphantom{1^{1^1}}}$
5004 \ & $4.0\cdot10^{13}$ \\
$10^{5}$ & $2.7\cdot10^{11}$ \\
$2\cdot10^{5}$ & $1.7\cdot10^{11}$ \\
$5\cdot10^{5}$ & $9.7\cdot10^{10}$ \\
$10^{6}$ & $4.4\cdot10^{10}$ \\
$2\cdot10^{6}$ & $2.6\cdot10^{10}$ \\
$8\cdot10^{6}$ & \ $4\cdot10^{9}$ \\
\hline
\end{tabular}
\end{table}

\newpage

\section{Missing: a prime}

By assuming $\rho$ is large when $p$ is small (a simplification!), we might ask for the least $p$, or rather,
``extremely large least $p$'' for given $q$ by searching for the first prime $p = kq+r$ in an AP for each residue
class $r$ coprime to $q$.  Li, Pratt, and Shakan \cite{LiPrattShakan2017} recently studied these primes for $q \leq 10^6$.
We doubled down on these efforts and went up to $q = 8.8\cdot10^6$ for even $q$,\footnote{\,Status at the time of writing,
the computation is still progressing indefinitely.} and found some especially large
first primes listed in table \ref{table7}. For comparison, a (semi-)conventional measure $\rho^* = p/(\varphi(q)\log^2 p)$
is employed here, while Li, Pratt, and Shakan use $p / (\varphi(q) \log \varphi(q) \log q)$ in their work.
The latter expression gives different results for odd $q$ when the same gap is achieved with $2q$ while the former doesn't. \\

As against the prime gaps (whether in AP or not), a lower bounding prime $p^-$ is missing here. To avoid any sort of
complication in using $\rho$ for small values of $r$ (including the fact that $G'(1)$ is undefined as \eqref{eq1} has a
simple pole at this point), in this case we substitute $p$ for $r$ in equation \eqref{eq2}, and further adapt equation
\eqref{eq3} with $p$ instead of $g$ in the denominator. \\

So we see that large $\rho$ are not necessarily easier to find when $p$, i.e. technically $p^+$, is as small as possible
(after all, we only look for the very first gap for each $r$, instead of all the other gaps with larger $p^-$), an
insight which also hints to a possible bound on $\rho$.

\begin{table}[H]
\caption{Data for least primes $p$ in AP which are largest for given $q$. Only the first exceptionally large ``gap'' and
subsequent champion ratios are shown, including all known examples where $\rho^* > 1$, meaning $p > \varphi(q)\log^2 p$.
\label{table7}
}
\centering
\small{
\begin{tabular}{ccccccl}
\hline 
\ $p$& $q {\vphantom{\displaystyle\int}}$ & $r$ & $k$ & $\rho$ &\!\!\!${p/(\varphi(q)\log^2 p)}$ & {\it Remark} \\
\hline
${\vphantom{1^{1^1}}}$
3 & 2 & 1 & 1 & 0.6267 & 2.4856 & skewed $\rho^*$ \\
2183963 & 23636 & 9451 & 92 & 1.0499 & 0.9185 & $1^{st}$ instance $\rho > 1$ \\
15714509 & 183336 & 130949 & 85 & 1.0406 & 0.9367 & new record $\rho^*$ \\
27361751 & 199432 & 39567 & 137 & 1.0706 & 0.9491 & new record $\rho$ and $\rho^*$ \\
136749709 & 783968 & 339277 & 174 & 1.1063 & 0.9941 & new record $\rho$ and $\rho^*$ \\
121770989 & 1084632 & 292205 & 112 & 1.1080 & 0.9958 & new record $\rho$ and $\rho^*$ \\
281309257 & 2732760 & 2567737 & 102 & 1.1101 & 1.0200 & $1^{st}$ instance $\rho^* > 1$ for $q > 2$ \\
673415261 & 4871052 & 1210085 & 138 & 1.1143 & 1.0094 & new record $\rho$ \\
1134128197 & 5497388 & 1666269 & 206 & 1.1508 & 1.0437 & new record $\rho$ and $\rho^*$ \\
1638279983 & 7197572 & 4431139 & 227 & 1.1096 & 1.0113 & $4^{th}$ of 5 known $\rho^* > 1$ for $q > 2$ \\
1653931571 & 8620136 & 7485595 & 191 & 1.1041 & 1.0112 & $5^{th}$ of 5 known $\rho^* > 1$ for $q > 2$ \\
\hline
\end{tabular}
}
\end{table}

\newpage

\section{Epilogue: data vs. proof}

Although we still seem to be far from proving it, the data corroborates that the gaps between
primes below $x$ are bounded by a constant times $\log^2x$. Until number theorists are able to find a search strategy
dramatically different from current methods and algorithms, it appears we cannot tell for sure whether or not
there are gaps between consecutive primes exceeding $\log^2p$, or more generally, whether or not there are infinitely
many gaps between primes in an arithmetic progression exceeding $\varphi(q) \log^2p$ for a given $q$.

\vspace{8mm}

\section{Acknowledgements}
Many thanks to Alexei Kourbatov for the initial inspiration and a lot of helpful exchange. The author would also like to thank
Andrew Granville for comments on a draft of the paper, and Bobby Jacobs from {\tt Mersenneforum.org} for moral support.

\vspace{20mm}

{
{\bf About the author:}
Martin Raab works in computational number theory as a hobby. He performed various computations regarding
prime numbers, some of which are found in the On-Line Encyclopedia of Integer Sequences (\url {https://oeis.org/}),
at the Mersenneforum (\url {https://www.mersenneforum.org/index.php}), and in the tables of first occurrence prime
gaps (currently hosted at \url {https://primegap-list-project.github.io/}). Martin lives in 63840 Hausen, Germany.
He can be reached at kilroy14159265@gmail.com.
}

\end{document}